# A Comprehensive Study of Complete Generalized New Mock Theta Functions

[1]Swayamprabha Tiwari, [2]Dr. Sameena Saba

Department of Mathematics & Statistics, Faculty of Science,

Integral University, Lucknow-226026

[1]E-mail:swayam6288@gmail.com

[2]E-mail:sameena9554@gmail.com

## Abstract

The generalization of new mock theta functions of Andrew's and Bringmann-et-al are given. Further we have given the expansion of these bilateral generalized new mock theta functions as $_2\phi_1$ series by Slater's transformation. After that we have given the continued fraction representation of these generalized mock theta functions.

## 1) Introduction

In his last letter to G .H. Hardy [8], S. Ramanujan listed seventeen mock theta functions of order three, five, six and seven. According to Ramanujan a mock theta function is a function $f(q), |q| < 1$, satisying the following two conditions:

(0) For every root of unity $\xi$, there is a $\theta$-function $\theta_\xi(q)$ such that the difference $f(q) - \theta_\xi(q)$ is bounded as $q \to \xi$ radially.

(1) There is no single $\theta$-function which works for all $\xi$, i.e., for every $\theta$-function $\theta(q)$ there is some root of unity $\xi$ for which $f(q) - \theta(q)$ is unbounded as $q \to \xi$ radially.

After studying mock theta functions, Andrews [3] and Bringmann, Hikami and Lovejoy [12] generated some new mock theta functions. In which four functions of Andrews and two of Bringmann, Hikami and Lovejoy found interesting. We have made comprehensive study of these mock theta functions [9-11].

The four mock theta functions of Andrews [3]:

$$\bar{\psi}_0(q) = \sum_{n=0}^{\infty} \frac{q^{2n^2}}{(-q;q)_{2n}} \qquad (1.1)$$

$$\bar{\psi}_1(q) = \sum_{n=0}^{\infty} \frac{q^{2n^2+2n}}{(-q;q)_{2n+1}} \qquad (1.2)$$

$$\bar{\psi}_2(q) = \sum_{n=0}^{\infty} \frac{q^{2n^2+2n}(q;q^2)_n}{(q^2;q^2)_n(-q;q)_{2n}} \qquad (1.3)$$

$$\bar{\psi}_3(q) = \sum_{n=0}^{\infty} \frac{q^{n^2}(-q;q)_n^2}{(q;q)_{2n}} \qquad (1.4)$$

The two mock theta functions of Bringmann , Hikami and Lovejoy[12] :

$$\bar{\phi}_0(q) = \sum_{n=0}^{\infty} q^n (-q;q)_{2n+1} \tag{1.5}$$

$$\bar{\phi}_1(q) = \sum_{n=0}^{\infty} q^n (-q;q)_{2n} \tag{1.6}$$

If we extend the definition of the defining series of a mock theta function to negative *n*, the series will become bilateral and we have Bilateral Mock Theta Functions. Watson called them "Complete".

## 2) Notations:

In this paper we have use the following notations:

$$(a;q^k)_n = \prod_{m=1}^{n}(1 - aq^{k(m-1)}), |q^k| < 1, \text{ n, non-negative integer}$$

$$(a;q^k)_0 = 1$$

$$(a;q^k)_\infty = \prod_{m=1}^{n}(1 - aq^{k(m-1)}),$$

$$(a_1, a_2, a_3 \ldots \ldots a_m; q^k)_n = (a_1;q^k)_n (a_2;q^k)_n \ldots \ldots (a_m;q^k)_n$$

$$_A\phi_{A-1}\begin{bmatrix} a_1, a_2, \ldots \ldots, a_A ; \\ b_1, b_2, \ldots \ldots, b_{A-1}; q_1, z \end{bmatrix}$$

$$= \sum_{n=0}^{\infty} \frac{(a_1;q_1)_n \ldots \ldots (a_A;q_1)_n \ z^n}{(b_1;q_1)_n \ldots \ldots (b_{A-1};q_1)_n \ (q_1;q_1)_n}, |z| < 1$$

## 3) Generalized New mock theta functions:

We have given the Generalization of new mock theta functions (1.1) – (1.6) as follows:

$$\bar{\psi}_0(t,\alpha,z;q) = \frac{1}{(t)_\infty} \sum_{n=0}^{\infty} \frac{(t)_n \ q^{2n^2-3n+n\alpha} \ z^{2n}}{\left(\frac{-z^2}{q};q\right)_{2n}} \tag{3.1}$$

$$\bar{\psi}_1(t,\alpha,z;q) = \frac{1}{(t)_\infty} \sum_{n=0}^{\infty} \frac{(t)_n \ q^{2n^2-n+n\alpha} \ z^{2n}}{\left(\frac{-z^2}{q};q\right)_{2n+1}} \tag{3.2}$$

$$\bar{\psi}_2(t,\alpha,z;q) = \frac{1}{(t)_\infty} \sum_{n=0}^{\infty} \frac{(t)_n \ q^{2n^2+2n-2n\alpha} \ z^{2n} \ (z;q^2)_n}{(z^2;q^2)_n \left(\frac{-z^2}{q};q\right)_{2n}} \tag{3.3}$$

$$\bar{\psi}_3(t,\alpha,z;q) = \frac{1}{(t)_\infty} \sum_{n=0}^{\infty} \frac{(t)_n \ q^{n^2-n+n\alpha} \ z^{2n} \ (-z;q)_n^2}{\left(\frac{-z^2}{q};q\right)_{2n}} \tag{3.4}$$

$$\bar{\phi}_0(t,\alpha,z;q) = \frac{1}{(t)_\infty} \sum_{n=0}^{\infty} (t)_n \ q^{n-2n\alpha} \ z^{2n} \left(\frac{-z^2}{q};q\right)_{2n+1} \tag{3.5}$$

$$\bar{\phi}_1(t,\alpha,z;q) = \frac{1}{(t)_\infty} \sum_{n=0}^{\infty} (t)_n \ q^{n-2n\alpha} \ z^{2n} \left(\frac{-z^2}{q};q\right)_{2n} \tag{3.6}$$

Now the complete form of new mock theta functions are

$$\bar{\psi}_{0,c}(t,\alpha,z;q) = \frac{1}{(t)_\infty} \sum_{n=-\infty}^{\infty} \frac{(t)_n q^{2n^2-3n+n\alpha} z^{2n}}{\left(\frac{-z^2}{q};q\right)_{2n}} \qquad (3.7)$$

$$\bar{\psi}_{1,c}(t,\alpha,z;q) = \frac{1}{(t)_\infty} \sum_{n=-\infty}^{\infty} \frac{(t)_n q^{2n^2-n+n\alpha} z^{2n}}{\left(\frac{-z^2}{q};q\right)_{2n+1}} \qquad (3.8)$$

$$\bar{\psi}_{2,c}(t,\alpha,z;q) = \frac{1}{(t)_\infty} \sum_{n=-\infty}^{\infty} \frac{(t)_n q^{2n^2+2n-2n\alpha} z^{2n} (z;q^2)_n}{(z^2;q^2)_n \left(\frac{-z^2}{q};q\right)_{2n}} \qquad (3.9)$$

$$\bar{\psi}_{3,c}(t,\alpha,z;q) = \frac{1}{(t)_\infty} \sum_{n=-\infty}^{\infty} \frac{(t)_n q^{n^2-n+n\alpha} z^{2n} (-z;q)_n^2}{\left(\frac{-z^2}{q};q\right)_{2n}} \qquad (3.10)$$

$$\bar{\phi}_{0,c}(t,\alpha,z;q) = \frac{1}{(t)_\infty} \sum_{n=-\infty}^{\infty} (t)_n\ q^{n-2n\alpha} z^{2n} \left(\frac{-z^2}{q};q\right)_{2n+1} \qquad (3.11)$$

$$\bar{\phi}_{1,c}(t,\alpha,z;q) = \frac{1}{(t)_\infty} \sum_{n=-\infty}^{\infty} (t)_n\ q^{n-2n\alpha} z^{2n} \left(\frac{-z^2}{q};q\right)_{2n} \qquad (3.12)$$

## 4) Expansion of Bilateral New Generalized Mock Theta Functions:

In the expansion of Slater [6, (5.4.3), p.129] taking $r=2$, we get:

$$\frac{\left(b_1,\ b_2,\ \frac{q}{a_1},\ \frac{q}{a_2},\ dz,\ \frac{q}{dz};q\right)_\infty}{\left(c_1,\ c_2,\ \frac{q}{c_1},\ \frac{q}{c_2};q\right)_\infty}\ {}_2\psi_2\left[\begin{matrix}a_1,\ a_2\\ b_1,\ b_2\end{matrix};q,z\right] =$$

$$\frac{q}{c_1}\frac{\left(\frac{c_1}{a_1},\frac{c_1}{a_2},\frac{qb_1}{c_1},\frac{qb_2}{c_1},\frac{dc_1z}{q},\frac{q^2}{dc_1z};q\right)_\infty}{\left(c_1,\frac{q}{c_1},\frac{c_1}{c_2},\frac{qc_2}{c_1};q\right)_\infty}\ {}_2\psi_2\left[\begin{matrix}\frac{qa_1}{c_1},\ \frac{qa_2}{c_1}\\ \frac{qb_1}{c_1},\ \frac{qb_2}{c_1}\end{matrix};q,z\right] + idem\ (c_1;c_2),\ (4.1)$$

Where $d = \frac{a_1 a_2}{c_1 c_2}$, $\left|\frac{b_1 b_2}{a_1 a_2}\right| < |z| < 1$

Expansion for equation $(3.7) - (3.12)$ bilateral new generalized mock theta functions are given as follows:

(i) Letting $t \to 0,\ q \to q^2, a_1 a_2 \to \infty, b_1 = \frac{-z^2}{q}, b_2 = -z^2, d = \frac{a_1 a_2}{c_1 c_2}, z \to \frac{z^2 q^{\alpha-1}}{a_1 a_2}$ in (4.1)

we get:

$$\frac{\left(\frac{-z^2}{q},-z^2,\frac{z^2 q^{\alpha-1}}{c_1 c_2},\frac{c_1 c_2}{z^2 q^{\alpha-3}};q^2\right)_\infty}{\left(c_1,c_2,\frac{q^2}{c_1},\frac{q^2}{c_2};q^2\right)_\infty} \cdot \bar{\psi}_{0,c}(0,\alpha,z;q)$$

$$= \frac{q^2}{c_1}\frac{\left(\frac{-z^2 q}{c_1},\frac{-z^2 q^2}{c_1},\frac{z^2 q^{\alpha-1}}{c_1 c_2},\frac{c_1 q^{5-\alpha}}{z^2};q^2\right)_\infty}{(c_1,\frac{q^2}{c_1},\frac{c_1}{c_2},\frac{q^2 c_2}{c_1};q^2)_\infty} \cdot \frac{1}{c_1^{2n}} \sum_{n=-\infty}^{\infty} \frac{q^{2n^2+n+n\alpha} z^{2n}}{\left(\frac{-z^2 q}{c_1};q^2\right)_n \left(\frac{-z^2 q^2}{c_1};q^2\right)_n} + idem\ (c_1,c_2)$$

$$(4.2)$$

(ii) Letting $t \to 0, q \to q^2, a_1 a_2 \to \infty, b_1 = -z^2, b_2 = -z^2 q, d = \frac{a_1 a_2}{c_1 c_2}, z \to \frac{z^2 q^{\alpha+1}}{a_1 a_2}$ in (4.1) we get:

$$\frac{\left(-z^2, -z^2 q, \frac{z^2 q^2}{c_1 c_2}, \frac{c_1 c_2}{z^2}; q^2\right)_\infty}{\left(c_1, c_2, \frac{q^2}{c_1}, \frac{q^2}{c_2}; q^2\right)_\infty} \cdot \bar{\psi}_{1,c}(0, \alpha, z; q)$$

$$= \frac{q^2}{c_1} \frac{\left(\frac{-z^2 q^2}{c_1}, \frac{-z^2 q^3}{c_1}, \frac{z^2}{c_2}, \frac{c_1}{z^2}; q^2\right)_\infty}{(c_1, \frac{q^2}{c_1}, \frac{c_1}{c_2}, \frac{q^2 c_2}{c_1}; q^2)_\infty} \cdot \frac{1}{c_1^{2n}} \sum_{n=-\infty}^{\infty} \frac{q^{2n^2+3n+n\alpha} z^{2n}}{\left(\frac{-z^2 q^2}{c_1}; q^2\right)_n \left(\frac{-z^2 q^3}{c_1}; q^2\right)_n} + \text{idem } (c_1, c_2)$$

(4.3)

(iii) Letting $t \to 0, q \to q^2, a_1 = -z^2, a_2 = -z^2 q, b_1 = b_2 = 0, d = \frac{a_1 a_2}{c_1 c_2}, z \to z^2 q^{1-2\alpha}$ in (4.1) we get:

$$\frac{\left(\frac{-q^2}{z^2}, \frac{-q}{z^2}, \frac{z^6 q^{2(1-\alpha)}}{c_1 c_2}, \frac{c_1 c_2}{z^4 q^{-2\alpha}}; q^2\right)_\infty}{\left(c_1, c_2, \frac{q^2}{c_1}, \frac{q^2}{c_2}; q^2\right)_\infty} \cdot \bar{\phi}_{0,c}(0, \alpha, z; q)$$

$$= \frac{q^2}{c_1} \frac{\left(\frac{-c_1}{z^2}, \frac{-c_1}{z^2 q}, \frac{z^6}{q^{2\alpha} c_2}, \frac{c_1 q^{2\alpha}}{z^2}; q^2\right)_\infty}{\left(c_1, \frac{q^2}{c_1}, \frac{c_1}{c_2}, \frac{q^2 c_2}{c_1}; q^2\right)_\infty} \cdot \sum_{n=-\infty}^{\infty} q^{n-2n\alpha} \left(\frac{-z^2 q^2}{c_1}; q^2\right)_n \left(\frac{-z^2 q^3}{c_1}; q^2\right)_n z^{2n} + \text{idem}(c_1, c_2)$$

(4.4)

(iv) Letting $t \to 0, q \to q^2, a_1 = \frac{-z^2}{q}, a_2 = -z^2, b_1 = b_2 = 0, d = \frac{a_1 a_2}{c_1 c_2}, z \to z^2 q^{1-2\alpha}$ in (4.1) we get:

$$\frac{\left(\frac{-q^3}{z}, \frac{-q^2}{z^2}, \frac{z^6}{c_1 c_2 q^{2\alpha}}, \frac{c_1 c_2 q^{2+2\alpha}}{z^6}; q^2\right)_\infty}{\left(c_1, c_2, \frac{q^2}{c_1}, \frac{q^2}{c_2}; q^2\right)_\infty} \cdot \bar{\phi}_{1,c}(0, \alpha, z; q) = \frac{q^2}{c_1} \frac{\left(\frac{-c_1}{z^2}, \frac{-q c_1}{z^2}, \frac{z^6 q^{-2-2\alpha}}{c_2}, \frac{c_2 q^{2+2\alpha}}{z^6}; q^2\right)_\infty}{\left(c_1, \frac{q^2}{c_1}, \frac{c_1}{c_2}, \frac{q^2 c_2}{c_1}; q^2\right)_\infty}$$

$$\cdot \sum_{n=-\infty}^{\infty} q^{n-2n\alpha} \left(\frac{-z^2 q}{c_1}; q^2\right)_n \left(\frac{-z^2 q^2}{c_1}; q^2\right)_n z^{2n} + \text{idem}(c_1, c_2)$$

(4.5)

## 5) Expansion for $\bar{\psi}_{0,c}(t, \alpha, z; q)$ and $\bar{\psi}_{1,c}(t, \alpha, z; q)$ as ${}_2\phi_1$ Series:

Taking $j = 1,2$ in [6, (5.4.5).p.130] we get:

$$\frac{\left(\frac{q}{a_1}, \frac{q}{a_2}, dz, \frac{q}{dz}; q\right)_\infty}{\left(\frac{q}{b_1}, \frac{q}{b_2}; q\right)_\infty} {}_2\psi_2 \begin{bmatrix} a_1, a_2 \\ b_1, b_2 \end{bmatrix}; q, z =$$

$$\frac{q}{b_1} \frac{\left(q, \frac{b_1}{a_1}, \frac{b_1}{a_2}, \frac{db_1 z}{q}, \frac{q^2}{db_1 z}; q\right)_\infty}{\left(b_1, \frac{q}{b_1}, \frac{b_1}{b_2}; q\right)_\infty} {}_2\phi_1 \begin{bmatrix} \frac{qa_1}{b_1}, \frac{qa_2}{b_1} \\ \frac{qb_2}{b_1} \end{bmatrix}; q, z + \text{idem } (b_1; b_2)$$

Where $d = \frac{a_1 a_2}{b_1 b_2}$. (5.1)

By using (5.1) we have given expansion of bilateral new generalized mock theta functions $(3.7) - (3.8)$ as $_2\phi_1$ series:

(i) Letting $q \to q^2$, $a_1 a_2 \to \infty$, $t \to 0$, $b_1 = \frac{-z^2}{q}$, $b_2 = -z^2$ and $z \to \frac{z^2 q^{\alpha-1}}{a_1 a_2}$ in (5.1) we get:

$$\bar{\psi}_{0,c}(t,\alpha,z;q) = \frac{-q^3}{z^2} \frac{\left(q^2, \frac{-z^4 q^{\alpha-4}}{c_1 c_2}, \frac{-q^{6-\alpha} c_1 c_2}{z^4}, \frac{-q^2}{z^2};q^2\right)_\infty}{\left(\frac{-z^2}{q}, \frac{1}{q}, \frac{z^2 q^{\alpha-1}}{c_1 c_2}, \frac{c_1 c_2 q^{3-\alpha}}{z^2};q^2\right)_\infty}$$

$$\times \sum_{n=0}^{\infty} \frac{q^{2n^2+3n+n\alpha} z^{-2n}}{(q^3;q^2)_n (q^2;q^2)_n} + \text{idem}(b_1,b_2) \qquad (5.2)$$

(ii) Letting $q \to q^2$, $a_1 a_2 \to \infty$, $t \to 0$, $b_1 = -z^2$, $b_2 = -z^2 q$ and $z \to \frac{z^2 q^{\alpha+1}}{a_1 a_2}$ in (5.1) we get:

$$\bar{\psi}_{1,c}(t,\alpha,z;q) = \frac{-q^2}{z^2} \frac{\left(q^4, \frac{-z^4 q^\alpha}{c_1 c_2}, \frac{-q^{3-\alpha} c_1 c_2}{z^4}, \frac{-q}{z^2};q^2\right)_\infty}{\left(-z^2, \frac{1}{q}, \frac{z^2 q^{\alpha+1}}{c_1 c_2}, \frac{c_1 c_2 q^{1-\alpha}}{z^2};q^2\right)_\infty} \times \sum_{n=0}^{\infty} \frac{q^{2n^2+3n+n\alpha} z^{-2n}}{(q^3;q^2)_n} + \text{idem}(b_1,b_2)$$

(5.3)

## 6) Expansion of $\bar{\phi}_{0,c}(t,\alpha,z;q)$ and $\bar{\phi}_{1,c}(t,\alpha,z;q)$ as $_2\phi_1$ Series:

For getting the expansion of a bilateral new generalized mock theta functions $(3.11) - (3.12)$ as a $_2\phi_1$ series, we take $c_j = qa_j, j = 1,2$ in Slater's expansion formula [6, (5.4.4), p.130] we get:

$$\frac{\left(b_1, b_2, \frac{q}{a_1}, \frac{q}{a_2}, z, \frac{q}{z};q\right)_\infty}{\left(qa_1, qa_2, \frac{1}{a_1}, \frac{1}{a_2};q\right)_\infty} \; _2\psi_2\begin{bmatrix} a_1, & a_2 \\ b_1, & b_2 \end{bmatrix};q,z]$$

$$= \frac{a_1\left(q, \frac{qa_1}{a_2}, \frac{b_1}{a_1}, \frac{b_2}{a_1}, a_1 z, \frac{q}{a_1 z};q\right)_\infty}{\left(qa_1, \frac{1}{a_1}, \frac{a_1}{a_2}, \frac{qa_2}{a_1};q\right)_\infty} \; _2\phi_1\begin{bmatrix} \frac{qa_1}{b_1}, & \frac{qa_1}{b_2} \\ \frac{qa_1}{a_2} \end{bmatrix}; q, \frac{b_1 b_2}{a_1 a_2 z}] + \text{idem}(a_1,a_2), \qquad (6.1)$$

where $\left|\frac{b_1 b_2}{a_1 a_2}\right| < |z| < 1$

(i) Letting $q \to q^2, b_1 = b_2 = 0, a_1 = -z^2, a_2 = -z^2 q, d = \frac{a_1 a_2}{c_1 c_2}, z \to z^2 q^{1-2\alpha}$ in (6.1), we have

$$\bar{\phi}_{0,c}(t,\alpha,z;q) = \frac{-z^2\left(q^2, q, -z^4 q^{1-2\alpha}, \frac{-q^{1+2\alpha}}{z^4}, -z^2 q^3, \frac{-1}{z^2 q};q^2\right)_\infty}{\left(\frac{1}{q}, q^3, \frac{-q^2}{z^2}, \frac{-q}{z^2}, z^2 q^{1-2\alpha}, \frac{q^{1+2\alpha}}{z^2};q^2\right)_\infty}$$

$$\times \sum_{n=0}^{\infty} \frac{q^{n-2n\alpha} z^{2n}}{(q;q^2)_n} + \text{idem}(a_1,a_2) \qquad (6.2)$$

(ii) Letting $q \to q^2, b_1 = b_2 = 0, a_1 = \frac{-z^2}{q}, a_2 = -z^2, d = \frac{a_1 a_2}{c_1 c_2}, z \to z^2 q^{1-2\alpha}$ in (6.1), we have

$$\bar{\phi}_{1,c}(t,\alpha,z;q) = \frac{-z^2}{q} \frac{\left(q^2, q, -z^2 q^{-2\alpha}, \frac{-q^{2+2\alpha}}{z^4}, -z^2 q^2, \frac{-1}{z^2};q^2\right)_\infty}{\left(\frac{1}{q}, q^3, \frac{-q^2}{z^2}, \frac{-q}{z^2}, z^2 q^{1-2\alpha}, \frac{q^{1+2\alpha}}{z^2};q^2\right)_\infty}$$

$$\times \sum_{n=0}^{\infty} \frac{q^{n-2n\alpha} z^{2n}}{(q;q^2)_n} + \text{idem}(a_1, a_2) \qquad (6.3)$$

## 7) Bilateral Generalized New Mock Theta Functions as continued fractions:

By using the continued fraction [1, (3.79), p.82] given by

$$\frac{\sum_{n=0}^{\infty} \frac{q^{n^2} \lambda^n}{(q)_n (-\beta)_n}}{\sum_{n=0}^{\infty} \frac{q^{n^2+n} \lambda^n}{(q)_n (-\beta)_n}} = 1 + \frac{\lambda q + \beta}{1-\beta} \cdot \frac{\lambda q^2 + \beta}{1-\beta \ldots}, \qquad (7.1)$$

We have given the continued fraction representation for $\bar{\psi}_{0,c}(t,\alpha,z;q)$, $\bar{\psi}_{1,c}(t,\alpha,z;q)$.

(i) Representation of $\bar{\psi}_{0,c}(t,\alpha,z;q)$ as a continued fraction.

Now (5.2) gives

$$\bar{\psi}_{0,c}(t,\alpha,z;q) = \frac{-q^3}{z^2} \frac{\left(q^2, \frac{-z^4 q^{\alpha-4}}{c_1 c_2}, \frac{-q^{6-\alpha} c_1 c_2}{z^4}, \frac{-q^2}{z^2} ; q^2\right)_\infty}{\left(\frac{-z^2}{q}, \frac{1}{q}, \frac{z^2 q^{\alpha-1}}{c_1 c_2}, \frac{c_1 c_2 q^{3-\alpha}}{z^2} ; q^2\right)_\infty} \sum_{n=0}^{\infty} \frac{q^{2n^2+3n+n\alpha}}{(q^3;q^2)_n (q^2;q^2)_n}$$

$$+ \frac{-q^2}{z^2} \frac{\left(q^2, \frac{-z^2 q^{\alpha-2}}{c_1 c_2}, \frac{-q^{5-\alpha} c_1 c_2}{z^4}, \frac{-q^3}{z^2} ; q^2\right)_\infty}{\left(z^2, \frac{z^2 q^{\alpha-1}}{c_1 c_2}, \frac{c_1 c_2 q^{3-\alpha}}{z^2}, q ; q^2\right)_\infty} \sum_{n=0}^{\infty} \frac{q^{2n^2+n+n\alpha}}{(q;q^2)_n (q^2;q^2)_n}$$

$$= S \sum_{n=0}^{\infty} \frac{q^{2n^2+3n+n\alpha}}{(q^3;q^2)_n (q^2;q^2)_n} + T \sum_{n=0}^{\infty} \frac{q^{2n^2+n+n\alpha}}{(q;q^2)_n (q^2;q^2)_n} \qquad (7.2)$$

Where

$$S = \frac{-q^3}{z^2} \frac{\left(q^2, \frac{-z^4 q^{\alpha-4}}{c_1 c_2}, \frac{-q^{6-\alpha} c_1 c_2}{z^4}, \frac{-q^2}{z^2} ; q^2\right)_\infty}{\left(\frac{-z^2}{q}, \frac{1}{q}, \frac{z^2 q^{\alpha-1}}{c_1 c_2}, \frac{c_1 c_2 q^{3-\alpha}}{z^2} ; q^2\right)_\infty}$$

and

$$T = \frac{-q^2}{z^2} \frac{\left(q^2, \frac{-z^2 q^{\alpha-2}}{c_1 c_2}, \frac{-q^{5-\alpha} c_1 c_2}{z^4}, \frac{-q^3}{z^2} ; q^2\right)_\infty}{\left(z^2, \frac{z^2 q^{\alpha-1}}{c_1 c_2}, \frac{c_1 c_2 q^{3-\alpha}}{z^2}, q ; q^2\right)_\infty}$$

Dividing the equation (7.2) by the first series on the right hand side of (7.2) we get;

$$\frac{\bar{\psi}_{0,c}(t,\alpha,z;q)}{\sum_{n=0}^{\infty} \frac{q^{2n^2+3n+n\alpha}}{(q^3;q^2)_n (q^2;q^2)_n}} = S + T \frac{\sum_{n=0}^{\infty} \frac{q^{2n^2+n+n\alpha}}{(q;q^2)_n (q^2;q^2)_n}}{\sum_{n=0}^{\infty} \frac{q^{2n^2+3n+n\alpha}}{(q^3;q^2)_n (q^2;q^2)_n}} \qquad (7.3)$$

Taking $q \to q^2, \lambda \to q^{1+\alpha}, \beta = -q$ in (7.1), and then putting it in the quotient of the summation on the right side of (7.3) to get the continued fraction representation for $\bar{\psi}_{0,c}(t,\alpha,z;q)$ as:

$$\frac{(1-q)\,\bar{\psi}_{0,c}(t,\alpha,z;q)}{\sum_{n=0}^{\infty} \frac{q^{2n^2+3n+n\alpha}}{(q^3;q^2)_n (q^2;q^2)_n}} = S + T\left[1 + \frac{q^2-q}{(1+q)} \cdot \frac{q^4-q}{(1+q)+\cdots\ldots}\right] \qquad (7.4)$$

(i) Representation of $\bar{\psi}_{1,c}(t,\alpha,z;q)$ as a continued fraction similar as earlier one.

$$\bar{\psi}_{1,c}(t,\alpha,z;q) = \frac{-q^2}{z^2} \frac{\left(q^4, \frac{-z^4 q^\alpha}{c_1 c_2}, \frac{-q^{3-\alpha} c_1 c_2}{z^4}, \frac{-q}{z^2}; q^2\right)_\infty}{\left(-z^2, \frac{1}{q}, \frac{-z^2 q^{\alpha+1}}{c_1 c_2}, \frac{c_1 c_2 q^{1-\alpha}}{z^2}; q^2\right)_\infty} \sum_{n=0}^\infty \frac{q^{2n^2+3n+n\alpha}}{(q^3;q^2)_n (q^2;q^2)_n}$$

$$+ \frac{-q}{z^2} \frac{\left(q^2, \frac{-z^4 q^\alpha}{c_1 c_2}, \frac{-q^{2-\alpha} c_1 c_2}{z^4}, \frac{-q^2}{z^2}; q^2\right)_\infty}{\left(z^2 q, q, \frac{z^2 q^{\alpha+1}}{c_1 c_2}, \frac{c_1 c_2 q^{1-\alpha}}{z^2}; q^2\right)_\infty} \sum_{n=0}^\infty \frac{q^{2n^2+n+n\alpha}}{(q;q^2)_n (q^2;q^2)_n}$$

$$\bar{\psi}_{1,c}(t,\alpha,z;q) = S_1 \sum_{n=0}^\infty \frac{q^{2n^2+3n+n\alpha}}{(q^3;q^2)_n (q^2;q^2)_n} + T_1 \sum_{n=0}^\infty \frac{q^{2n^2+n+n\alpha}}{(q;q^2)_n (q^2;q^2)_n}$$

Where

$$S_1 = \frac{-q^2}{z^2} \frac{\left(q^4, \frac{-z^4 q^\alpha}{c_1 c_2}, \frac{-q^{3-\alpha} c_1 c_2}{z^4}, \frac{-q}{z^2}; q^2\right)_\infty}{\left(-z^2, \frac{1}{q}, \frac{-z^2 q^{\alpha+1}}{c_1 c_2}, \frac{c_1 c_2 q^{1-\alpha}}{z^2}; q^2\right)_\infty}$$

$$T_1 = \frac{-q}{z^2} \frac{\left(q^2, \frac{-z^4 q^\alpha}{c_1 c_2}, \frac{-q^{2-\alpha} c_1 c_2}{z^4}, \frac{-q^2}{z^2}; q^2\right)_\infty}{\left(z^2 q, q, \frac{z^2 q^{\alpha+1}}{c_1 c_2}, \frac{c_1 c_2 q^{1-\alpha}}{z^2}; q^2\right)_\infty}$$

Hence we have

$$\frac{(1+q)\, \bar{\psi}_{1,c}(t,\alpha,z;q)}{\sum_{n=0}^\infty \frac{q^{2n^2+3n+n\alpha}}{(q^3;q^2)_n (q^2;q^2)_n}} = T_1 + S_1 \left[1 + \frac{q^{2+\alpha}-q}{1+q} \cdot \frac{q^{3+\alpha}-q}{1+q}\right] \qquad (7.5)$$

## Conclusion

The study of bilateral generalized new mock theat functions is interesting as it gives the expansion by Slater and also expansion for these generalized new mock theta functions as $_2\phi_1$ Series. After that we have given continued fractions for generalized new mock theta functions.

## Acknowledgement

I am thankful to my guide in my research work . My MCN is IU/R&D/2023-MCN0002042.